\documentclass{amsart}
\usepackage{amsmath, amssymb, xy}
\input xy
\xyoption{all}
\usepackage{pdflscape}
\usepackage{hyperref}

\numberwithin{equation}{section}

\newtheorem{thm}{Theorem}[section]
\newtheorem{prop}[thm]{Proposition}
\newtheorem{conj}[thm]{Conjecture}
\newtheorem{cor}[thm]{Corollary}
\newtheorem{lem}[thm]{Lemma}

\theoremstyle{definition}
\newtheorem{defn}[thm]{Definition}
\newtheorem*{remark*}{Remark}
\newtheorem{remark}[thm]{Remark}
\newtheorem*{pf}{Proof}

\newcommand{\F}{\mathbf{F}}

\newcommand{\cO}{\mathcal{O}}

\newcommand{\cA}{\mathcal{A}}
\newcommand{\D}{\mathcal{D}}
\newcommand{\fp}{\kappa}
\newcommand{\s}{\mathfrak{s}}

\newcommand{\OO}[1]{\cO_{#1}}
\newcommand{\G}{\mathrm{G}}
\newcommand{\GL}{\text{GL}}
\newcommand{\SL}{\text{SL}}

\newcommand{\glO}[2]{\text{GL}_{#1}({\cO}_{#2})}
\newcommand{\mO}[2]{\text{M}_{#1}({\cO}_{#2})}

\newcommand{\ind}{\text{Ind}}

\newcommand{\res}{\text{Res}}
\newcommand{\tr}{\text{Tr}}
\newcommand{\irr}{\mathrm{Irr}}
\title[On representations of General linear groups]{On
  representations of General linear groups over
  principal ideal local rings of length two }
\author{Pooja Singla}
\address{Institute of Mathematical Sciences, Chennai, India 600113}
\email{pooja@imsc.res.in}
\begin{document}
\begin{abstract}
  We study the irreducible complex representations of general linear groups over principal ideal local rings of length two with a fixed finite residue field.
  We construct a canonical correspondence between the irreducible representations of all such groups which preserves dimensions.
  For general linear groups of order three and four over these rings, we construct all the irreducible representations.
  We show that the the problem of constructing all the irreducible representations of all general linear groups over these rings is not easier than the problem of constructing all the irreducible representations of the general linear groups over principal ideal local rings of arbitrary length in the function field case.
\end{abstract}
\maketitle
\tableofcontents
\section{Introduction}
Let $\mathrm{F}$ be a non-Archimedean local field with ring of integers $\cO$.
Let $\wp$ be the unique maximal ideal of $\cO$ and $\pi$ be a fixed
uniformizer of $\wp$.
Assume that the residue field $\cO/\wp$ is finite.
The typical examples of such rings of integers are
$\mathbf Z_{p}$ (the ring of p-adic integers) and $\mathbf F_q [[t]]$ (the ring of
formal power series with coefficients over a finite field). We denote by $\cO_{\ell}$
the reduction of $\cO$ modulo $\wp^{\ell}$, i.e. $\cO_{\ell} =
\cO/\wp^{\ell}$. Let $\Lambda_{k}$ denote the set of partitions
with $k$ parts, namely, non-increasing finite sequences $(\ell_1,
\ell_2, \ldots, \ell_k)$ of positive integers, and let $\Lambda = \cup  \Lambda_{k}$. Since
$\cO$ is a principal ideal domain with a unique maximal ideal $\wp$, every
finite $\cO$-module is of the form $\oplus_{i=1}^{k}\cO_{\ell_{i}}$, where
$\ell_{i}$'s can be arranged so that $\lambda = (\ell_1, \ell_2,
\ldots, \ell_k) \in
\Lambda_{k}$. Let $\mathrm{M}_{\lambda} = \oplus_{i=1}^{k} \cO_{\ell_{i}}$ and
\[
\G_{\lambda, F} = \mathrm{Aut}_{\mathcal{O}}(\mathrm{M}_{\lambda}).
\]
We write $G_{\lambda}$ instead of $G_{\lambda, F}$ whenever field $F$ is clear from the context. 
If $\mathrm{M}_{\lambda} = \cO_{\ell}^{n}$ for some natural number $n$, then
the group $\G_{\lambda}$ consists of invertible matrices of order $n$ with
entries in the ring $\cO_{\ell}$, so we use the notation $\glO{n}{\ell}$ for
$\G_{\lambda}$ in this case.

The representation theory of the finite groups
$G_{\lambda}$ has attracted the attention of many mathematicians. We give a brief 
history of this problem. Green~\cite{MR0072878}
calculated the characters of the irreducible representations of
 $\GL_{n}(\cO_1)$. Several authors, for instance,    
Frobenius~\cite {Frobenius2},  Rohrbach~\cite {Rohrbach}, 
 Kloosterman~\cite {MR0021032, MR0021033}, Tanaka~\cite {MR0229737},
 Nobs-Wolfart~\cite {MR0429742}, Nobs~\cite{MR0447489}, Kutzko~\cite {MR0320170},
Nagornyi~\cite{MR0491919} and Stasinski~\cite{stasinski-2008} studied the representations of
the groups  $\SL_{2}(\cO_{\ell})$ and $\glO{2}{\ell}$. Nagornyi~\cite{MR0498881} obtained 
partial results regarding the representations of $\glO{3}{\ell}$
and Onn~\cite{MR2456275} constructed all the irreducible
representations of the groups $\G_{(\ell_1, \ell_2)}$. Recently, Avni-Klopsch-Onn-Voll~\cite{avni-2009} have announced results about the representation theory of the groups $\mathrm{SL}_{3}(\mathbf Z_p)$.

In another direction, it was observed that, being maximal compact subgroups, $\GL_n(\cO)$ 
play an important role in the representation theory of the groups $\GL_n(F)$. Further, every continuous representation 
of $\GL_n(\cO)$ factors through one of the natural homomorphisms
$\GL_n(\cO) \rightarrow \GL_n(\cO_{\ell})$. This 
brings the study of irreducible representations of groups $\GL_n(\cO_{\ell})$ to the forefront.
Various questions regarding the complexity of the problem of determining irreducible representations of these groups 
were asked. For example Nagornyi~\cite{MR0498881} proved that 
that this problem contains matrix pair problem. Aubert-Onn-Prasad-Stasinski~\cite{AOPS} proved that, 
for $F = \F_q((t))$, constructing all irreducible representations of $\glO{n}{2}$ for all $n$ is
equivalent to constructing all irreducible representations
of $G_{\lambda, \F_{q^m}((t))}$ for all $\lambda$ and $m$ (see also Section~\ref{complexity}). 

Motivated by Lusztig's work for finite 
groups of Lie type, Hill~\cite{MR1231615}  partitioned all the irreducible 
representations of groups $\GL_n(\cO_{\ell})$ into geometric conjugacy classes and reduced
the study of irreducible representations of $\GL_n(\cO_{\ell})$ to the study of its nilpotent characters.
In later publications ~\cite{MR1265002, MR1334228, MR1311772}, he succeeded in constructing many 
irreducible representations (namely strongly-semisimple, 
semisimple, regular etc.) for these groups. Following the techniques
used in the representation theory of groups $\glO{n}{1}$ and $\GL_n(F)$, various notions
like cuspidality and supercuspidality were introduced for
representations of $\GL_n(\cO)$ (for more on this see ~\cite{AOPS}), but 
the complete knowledge of irreducible 
representations of groups $\GL_n(\cO_{\ell})$ for $\ell \geq 2$ is still unknown. From the available results, it was observed 
that methods of constructing irreducible
representations of groups $\G_{\lambda}$ do not depend on the particular
ring of integers $\cO$ and depend only on the residue field. This led Onn to
conjecture \cite[Conjecture 1.2]{MR2456275} that
\begin{conj}
The isomorphism type of the
group algebra $\mathbb C[G_{\lambda}]$ depends only on $\lambda$ and $q = |\cO/\wp|$. 
\end{conj}
  
We discuss the method of constructing complex irreducible
representations of the groups $\GL_{n}(\cO_{2})$ with the help of Clifford theory and reduce this problem 
to constructing irreducible representations of 
certain subgroups of $\glO{n}{1}$. This enables us to give an affirmative answer 
to the above conjecture for $\glO{n}{2}$. 
The groups $\glO{n}{2}$, for distinct rings of integers 
$\mathcal{O}$ are not necessarily isomorphic, even when the residue fields are isomorphic. 
For example; for a natural number $n$ and a prime
$p$, the group $\GL_{n}(\mathbf F_p[[t]]/t^{2})$ is a semi-direct
product of the groups $\mathrm{M}_{n}(\F_{p})$ and $\GL_{n}(\mathbf F_{p})$, but on the other hand
$\GL_{n}(\mathbf Z_p/p^2 \mathbf Z_p)$ is not unless $n = 1$ or $(n,p) = (2,2), (2,3)
\;\mathrm{or}\;(3,2)$
(Sah~\cite[p.~22]{MR0463264}, Ginosar~\cite{MR1872829}).
Our main emphasis is on proving that all of 
their irreducible representations can be constructed in a uniform way. We also succeed in 
showing that representation theory of groups $G_{\lambda, \mathbf F_{q^m}((t))}$ plays a 
vital role in  representation theory of groups
$\glO{n}{2}$ for any $\cO$, in the sense that if we know irreducible
representations of the groups $G_{\lambda, \mathbf F_{p^m}((t))}$ for all
positive integers $m$, we can
determine all the representations of $\glO{n}{2}$. 

More precisely,
let $\mathrm{F}$ and $\mathrm{F}'$ be local fields with rings of integers
  $\mathcal{O}$ and $\mathcal{O'}$ respectively such that their residue fields
  are finite and isomorphic (with a fixed isomorphism). Let $\wp$ and $\wp'$ be the
  maximal ideals of $\mathcal{O}$ and $\mathcal{O}'$ respectively. As
  described earlier, $\cO_2$ and $\cO'_2$ denote the rings
  $\cO/\wp^2$ and $\cO'/\wp'^2$ respectively.
We prove

\begin{thm}[Main Theorem]
%\marginpar{thm:main}
\label{thm:main} 
There exists a canonical bijection between the irreducible representations of
  $\GL_{n}(\mathcal{O}_2)$ and those of 
  $\GL_{n}(\mathcal{O}'_2)$, which preserves dimensions. 

\end{thm}

\begin{defn}(Representation Zeta function)
Let $G$ be a group. The representation zeta function of $G$ is the function
\[
R_{G}({\D}) = \sum_{\rho \in \irr{G}} \D^{\dim{\rho}} \in \mathbb Z[\mathcal{D}]
\] 
\end{defn}
\noindent In view of the above definition, Theorem~\ref{thm:main} implies that  
\begin{cor}
\label{cor:representation zeta polynomial}
The representation zeta functions of $\glO{n}{2}$ and $\GL_{n}(\cO'_2)$
are equal.
\end{cor}
In other words, the representation zeta function depends on the ring only through the order of its residue field.

Concerning the complexity of the problem of constructing irreducible representations of groups $\glO{n}{2}$,
we obtain the following generalisation of \cite[Theorem 6.1]{AOPS}. 
\begin{thm} Let $\cO$ be the ring of integers of a non-Archimedean local field
  $F$, such that residue field has cardinality $q$. 
  Then the problem of constructing irreducible representations of the
  following groups are equivalent:
\begin{enumerate}
\item $\glO{n}{2}$ for all $n \in \mathbf N$. 
\item $G_{\lambda, E}$ for all partitions $\lambda$ and all unramified
  extensions $E$ of $\mathbf F_q((t))$.   
\end{enumerate}
\end{thm} 

We construct all the irreducible representations of
$\glO{2}{2}, \glO{3}{2}$ and $\glO{4}{2}$. As mentioned earlier, the representation theory of
$\glO{2}{2}$ is already known. Partial results regarding the
representations of $\glO{3}{2}$ has been obtained by Nagornyi \cite{MR0498881}
but the representation theory of $\glO{4}{2}$ seems completely novel. We find that

\begin{thm}
%\marginpar{thm:zeta polynomial}
\label{thm:zeta polynomial} The number and dimensions of irreducible representations of groups
  $\glO{3}{2}$ and $\glO{4}{2}$ are polynomials in $\mathbb{Q}[q]$. 
\end{thm}

This theorem proves the strong version of Onn's conjecture~\cite[Conjecture~1.3]{MR2456275} for the groups $\GL_{3}(\cO_2)$ and $\GL_{4}(\cO_2)$.\\
\vspace*{1mm}

\subsection{Organization of the article}
In Section~$2$, we set up the basic notation that we use throughout the
article and discuss the action of the group $\glO{n}{2}$ on the characters of its
normal subgroup $K = \rm{Ker}(\glO{n}{2} \mapsto \glO{n}{1})$. We also state
the results of Clifford theory, which we use later to prove Theorem~\ref{thm:main}. 

In Section~$3$  we briefly review the similarity classes of
$\mathrm M_{n}(\mathbf F_{q})$ and in Section~$4$, we 
discuss the centralizer algebras of matrices, namely the set of matrices that
commute with a given matrix. For any matrix $A \in
\rm{M}_n(\mathbf F_q)$ in Jordan canonical form, we describe
explicitly its centralizer in $\mathrm{M}_{n}(\mathbf F_q)$ and in $\GL_{n}(\mathbf F_q)$.

In Section~$5$, we present the 
proof of Theorem~\ref{thm:main}. For this, we prove that all the characters of the subgroup
$K$ can be extended to its stabilizer in $\glO{n}{2}$. 

Section~$6$ is devoted to applications of Theorem~\ref{thm:main}. We
express a relation between the representation zeta function of
$\glO{n}{2}$ and that of centralizers in $\glO{n}{1}$. Then we describe the representation zeta functions of
$\glO{2}{2}, \glO{3}{2}$, $G_{(2,1,1)}$, and $\glO{4}{2}$. In particular we prove
Theorem~\ref{thm:zeta polynomial}.

\subsection{Acknowledgements} I thank my advisor Amritanshu
Prasad for his invaluable guidance and support throughout this
project. It is a great pleasure to thank Uri Onn for many
fruitful discussions and, in particular, for
the proof of
Lemma~\ref{lem:G(2,1,1)}, which he provided through a personal communication. I
also thank  Alexander Stasinski for reading first draft of this
article and providing very useful feedback.

\section{Notations and Clifford Theory}
%\marginpar{section:clifford theory}
\label{section:clifford theeory}
In this section we set up the basic notation that we use throughout the
article. We state and apply the main results of Clifford theory to our case and state
Proposition~\ref{prop:main}, which is
an important step towards the proof of Theorem~\ref{thm:main}\\

In this article, by character we mean a one dimensional
representation, unless stated otherwise. For any group $G$, we denote
by $\rm{Irr}(G)$, the set of its irreducible representations, and for any
abelian group $A$, we denote by $\hat{A}$, the set of its characters.

 Let $\fp: \glO{n}{2} \rightarrow \glO{n}{1}$ be the natural quotient
map and $K = \rm{Ker}(\fp)$. Then $A \mapsto I + \pi A$ induces an isomorphism $\mathrm M_{n}(\mathcal{O}_1) \tilde{\to}
K$. Fix a non-trivial additive character $\psi : \mathcal{O}_{1} \rightarrow
\mathbf{C}^{*}$ and for any $A \in \mO{n}{1}$ define $\psi_{A}: K
\rightarrow \mathbf{C}^{*}$ by 
\[
\psi_{A}(I + \pi X) = \psi(\mathrm{Tr}(AX)).
\]
Then $A \mapsto \psi_{A}$ gives an isomorphism $\mathrm M_{n}(\mathcal{O}_1)
\tilde{\to} \hat{K}$. The group $\GL_{n}(\mathcal{O}_2)$ acts on $\mathrm M_{n}(\mathcal{O}_1)$ by conjugation
via its quotient $\GL_{n}(\mathcal{O}_1)$, and therefore on $\hat{K}$: for $\alpha
\in \glO{n}{2}$ and $\psi_A \in \hat{K}$, we have
%\marginpar{eq:orbit} 
\begin{eqnarray}
\label{eq:orbit}
\psi^{\alpha}_{A}(I + \pi X) & = & \psi_A(I + \pi \alpha X \alpha^{-1}) {}
                          \nonumber  \\ 
                          & = &  \psi(\tr(A \fp(\alpha) X \fp(\alpha)^{-1}))
                                                   {}  \nonumber \\
                           &   = & \psi(\tr(\fp(\alpha)^{-1} A \fp(\alpha) X ))
                          \nonumber {}\\
                           &  =  &  \psi_{\fp(\alpha)^{-1}A \fp(\alpha)} (I + \pi X) 
\end{eqnarray}
Thus the action of $\glO{n}{2}$ on the characters of $K$ transforms to its conjugation
(inverse) action on elements of $\mO{n}{1}$.

We shall use the following results of Clifford theory.

\begin{thm}
%\marginpar{thm:clifford}
\label{thm:clifford} Let $G$ be a finite group and $N$ be a normal subgroup. For any
  irreducible representation $\rho$ of $N$, let $T(\rho) =
  \{ g \in G |\,\, \rho^{g} = \rho \}$ denote the stabilizer of $\rho$. Then
  the following hold 
\begin{enumerate}

\item If $\pi$ is an irreducible representation of $G$ such that $\langle
  \pi|_{N}, \rho \rangle \neq 0$, then $\pi|_{N} = e (
  \oplus _{\rho \in \Omega} \rho )$ where $\Omega$ is an orbit of irreducible
  representations of $N$ under the action of $G$, and $e$
  is a positive integer. \\

\item Suppose that $\rho$ is an irreducible representation of $N$. Let
$$
A = \{ \theta \in \irr (T(\rho)) | \langle \res ^{T(\rho)}_{N} \theta, \rho \rangle
\neq 0 \} $$

$$ B = \{ \pi \in \irr G | \langle \res^{G}_{N} \pi, \rho \rangle
\neq 0 \} $$
Then $$ \theta \rightarrow \ind^{G}_{T(\rho)}(\theta)$$ is a bijection
of $A$ onto $B$. \\

\item Let $H$ be a subgroup of $G$ containing $N$, and suppose that $\rho$ is
an irreducible representation of $N$ which has an extension $\tilde{\rho}$ to
$H$ (i.e. $\tilde{\rho}|_{N} = \rho $). Then the representations $\chi \otimes \rho$
for $\chi \in \irr(H/N)$ are irreducible, distinct for distinct $\chi$, and \\
\[
\ind^{H}_{N}(\rho) = \oplus_{\chi \in \irr(H/N)} \chi \otimes
\tilde{\rho}.
\]   
\end{enumerate}    
\end{thm}

For proofs of the above, see for example, 6.2, 6.11, and
6.17 respectively in Isaacs~\cite{MR0460423}.  Applying the above results to
the group $
G = \GL_{n}(\cO_{2})$ and normal subgroup $N = K$, we see that the following proposition plays an important role in the proof of Theorem~\ref{thm:main}.
\begin{prop}
%\marginpar{prop:main}
\label{prop:main} For a given $A \in \mathrm M_{n}(\mathcal{O}_1)$, there exists a 
character $\chi$ of $T(\psi_A)$ such that $\chi|_{K} = \psi_{A}$ (such a 
character $\chi$ is called an extension of $\psi_A$). 
\end{prop}

%%%%%%%%%%%%%%%%%%%%%%%%%%%%%%%%%%%%%%%%%%%%%%%%%%%%%%%%%%%%%%%
%primary decomposition and Jordan canonical form
%
%
%
%%%%%%%%%%%%%%%%%%%%%%%%%%%%%%%%%%%%%%%%%%%%%%%%%%%%%%%%%%%%%%%

\section{Primary Decomposition and Jordan Canonical Form}
In this section we describe the primary decomposition of matrices
under the action of conjugation. We also discuss
Jordan canonical form for those matrices whose characteristic polynomials split
over $\mathbf F_q$.\\ 

\noindent Let $f$ be an irreducible polynomial with coefficients in $\mathbf F_q$.  
\begin{defn}(f-primary Matrix) A matrix with entries in $\mathbf F_q$ is $f$-primary if its
  characteristic polynomial is a power of $f$. 
\end{defn}

\noindent {\bf Notation:} If $A_{i}$'s for $1 \leq i \leq l$ are matrices, then we
denote by $\oplus_{i}A_{i}$ the block diagonal matrix
\[
\begin{bmatrix}
A_{1} & 0  & \cdots & 0  \\
0 & A_{2} & \cdots       & 0 \\
\vdots & \vdots & \ddots &     \vdots \\
0      & 0 & \cdots      &  A_{l} 
\end{bmatrix}
\] 
In the same spirit, if $\cA_{1}, \cA_{2}, \ldots, \cA_{l}$ are the sets of matrices
then $\oplus_{i=1}^{l}
\mathcal{A}_{i}$ denotes the set of  matrices $\{ \oplus_{i=1}^{l} A_{i}| \;\; A_{i} \in
\mathcal{A}_{i} \}$. 
The next two theorems are easy consequences of structure theorems for
$\mathbf F_q[t]$-modules. For proofs of these see
Bourbaki~\cite[A.VII.31]{MR1994218}.   
\begin{thm}
%\marginpar{thm:primary}
\label{thm:primary}(Primary Decomposition) Every matrix $A \in \mathrm M_{n}(\mathbf F_q)$ is similar to a
  matrix of the form 
\[
\oplus_{f} A_{f}. 
\]

\noindent where $A_{f}$ is an f-primary matrix, and the sum is over the
irreducible factors of the characteristic polynomial of $A$. Moreover,
for every f, the similarity class of $A_{f}$ is uniquely determined by
the similarity class of $A$. 
\end{thm}

\begin{defn}(Split Matrix) A matrix over $\mathcal{O}_1$ is called
  split if its characteristic polynomial splits over 
  $\mathbf F_q$.
\end{defn}

\begin{defn}(Elementary Jordan Blocks)  For a natural number $n$ and an
  element $a \in \mathbf F_q$, elementary Jordan block $J_{n}(a)$ is the matrix
\[
\begin{bmatrix}
a & 1 & 0 & 0 &  \cdots & 0 \\
0 & a & 1 &  0 &  \cdots & 0 \\
0 & 0 & a &  1&       &    \\
\vdots & \vdots & &  \ddots & \ddots & \vdots \\
       &        &        &  &    a  &        1 \\
0      &   0    &      \cdots  &   & 0 &   a 
\end{bmatrix}_{n \times n}
\]
\end{defn}

%\begin{defn}(Partition) Let $n$ be a natural number. A partition of $n$ is a tuple
%  $(a_{1}, a_{2}, ......, a_{k})$ of positive integers such that
%  $a_{1} \geq a_{2} \geq.....\geq a_{k}$ and 
%  $a_{1} + a_{2} + ...+ a_{k} = n$. 
%\end{defn}

\begin{thm}
%\marginpar{thm:jordan}
\label{thm:jordan} (Jordan Canonical Form for Split Matrices)  Every split matrix $A
  \in \mathrm M_{n}(\mathbf F_q)$, up to the rearrangement of the
$a_i$'s, is similar to a unique matrix of the form 
\[
\oplus_{i} J_{\lambda(a_i)}(a_{i})
\]
where $\lambda(a) = (\lambda_1(a), \lambda_2(a),\ldots, \lambda_k(a))$ is a
partition and 
\[
J_{\lambda(a_i)}(a_i) = 
\begin{bmatrix}
J_{\lambda_{1}(a_i)}(a_{i})  & 0 & \cdots & 0 \\
0             &  J_{\lambda_{2}(a_i)}(a_i) & \cdots & 0 \\
\vdots      &   \vdots                &  \ddots &  \vdots \\
0 &   0  &  \cdots  &  J_{\lambda_{k}(a_i)}(a_{i}) 
\end{bmatrix}
\]
and each $J_{\lambda_{j}(a_i)}(a_i)$ is an elementary Jordan block with
eigenvalue $a_{i}$.
\end{thm}

%%%%%%%%%%%%%%%%%%%%%%%%%%%%%%%%%%%%%%%%%%%%%%%%%%%%%%%%%%%%%%%%%
%%% Centralizers
%
%
%
%%%%%%%%%%%%%%%%%%%%%%%%%%%%%%%%%%%%%%%%%%%%%%%%%%%%%%%%%%%%%%%%

\section{Centralizers} Let $R$ be a commutative ring with unity. In this section we 
determine centralizers (see Def~\ref{def:centralizer}) of certain matrices in $\mathrm M_n(R)$ and $GL_n(R)$. 
We also relate the groups $G_{\lambda, \mathbf F_{q^m}((t))}$ with centralizers in $GL_n(\mathbf F_q)$.

\begin{defn}(Centralizer of an element)
%\marginpar{def:centralizer}
\label{def:centralizer}
 Let $L$ be a semigroup under
  multiplication and $l$ be an element of $L$. Assume that $T$ is a subset of $L$. Then
 { \em centralizer} of $l$ in $T$, $Z_{T}(l)$, is the set of elements of $T$ that
  commute with $l$, i.e.,
\[
Z_{T}(l) = \{ t \in T| \,\, tl = lt \}
\] 
\end{defn}
\begin{remark}
If $T$ is a group, then $Z_{T}(l)$ is a subgroup of $T$. 
\end{remark}
\begin{defn}(Principal Nilpotent Matrix) A square matrix is called Principal
  Nilpotent, if it is of the form
\[
N_{n} = 
\begin{bmatrix}
0 & 1 & 0 & \cdots & 0 \\
0 & 0 & 1 & 0   &    0 \\
\vdots &  & \ddots & \ddots & \vdots \\
0     &  \cdots   &    &  0  & 1\\
0    &   \cdots   &    &  0  & 0  
\end{bmatrix}_{n \times n}
\]
\end{defn}

In the sequel we use the notation $N_n$ for the principal Nilpotent matrix of
order $n$.\\

\vspace*{1mm}

\noindent  Let $n_1, n_2,\ldots, n_l$ be a
  sequence of natural numbers, such that $n =
  n_1 + n_2 + \cdots + n_l$. Let $A = \oplus_{i=1}^{l} N_{n_{i}}$. In 
Lemmas~\ref{lem:commuters}-\ref{lem:centralizer} we describe the centralizer algebras
  $Z_{\mathrm M_n(R)}(A)$. The proofs of these Lemmas involve simple
  matrix multiplications, so we leave these for reader.    
\begin{lem}
%\marginpar{lem:commuters}
\label{lem:commuters}
 Let $a$ and $b$ be two elements of $R$ such that $a-b$ is a unit in $R$. Assume that
  $A$ and $B$ are two upper triangular matrices such that all the
  diagonal entries of $A$ are equal to $a$ and those of $B$ are equal
  to $b$. Then there
  does not exist any non-zero matrix $X$ over $R$ such that $ X A = B
  X $. 
\end{lem}
\begin{lem}
%\marginpar{lem:generalcentralizers}
\label{lem:generalcentralizers}
 Let $a_{1}, a_{2},\ldots ,a_{l}$ be elements of $R$ such that for all $i \neq j$, $a_{i}- a_{j}$ is
 invertible in $R$. Let $A = \oplus_{i=1}^{l} A_{i}$ be a square
 matrix of order $n$, where $A_{i}$'s are upper triangular matrices of
 order $n_{i}$. Assume that all diagonal entries of  $A_{i}$ are equal
 to $a_{i}$. Then,  
\[
Z_{\GL_{n}(R)}(A) = \oplus_{i=1}^{l} Z_{\GL_{n_i}(R)} (A_{i}) 
\]
\end{lem}
\begin{defn}(Upper Toeplitz Matrix) A square matrix of order $n$ is called Upper Toeplitz, if it is of the form 
\[
\begin{bmatrix}
a_{1} & a_{2} & \cdots & a_{n-1} & a_{n} \\
0    & a_{1}  & a_{2}&    &    a_{n-1} \\
\vdots & \ddots   & \ddots   &  \ddots   &    \vdots \\
      &  \cdots   &       & a_{1}    &  a_{2} \\
0      &   \cdots   &  0    &   0      &   a_{1}   
\end{bmatrix}_{n \times n }
\]
\end{defn}
\begin{lem}
%\marginpar{lem:toeplitz} 
\label{lem:toeplitz}
Assume that $N_{n}$ and $N_{m}$ are principal nilpotent matrices of order $n$ and $m$
  respectively. Then the matrices $X$ over $R$ such that 
$ X N_{m} = N_{n} X $ are of the form
\[
X = \Bigg\{
\begin{array}{ccc}
\left[
\begin{array}{ll}
0_{n \times m-n} &  T_{n \times n}  
\end{array}
\right]
\textrm{ if $ n \leq m $} \\
              \\
\left[
\begin{array}{cc}
T_{m \times m} \\
0_{n-m \times m}
\end{array} 
\right]
\textrm{ if $ n \geq m $} 
\end{array}
\]
Where $T_{ s \times s }$, for a natural number $s$, is an upper Toeplitz Matrix
of order $s$, over the ring $R$. 
\end{lem}

This Lemma motivates the following definition of rectangular upper Toeplitz
matrix. 
\begin{defn}
%\marginpar{lem:rectangular}
\label{lem:rectangular}
(Rectangular Upper Toeplitz Matrix) A matrix of order $n
  \times m$, over a ring $R$ is called a
  rectangular upper Toeplitz if it is of the form
\[
\begin{array}{lll}
\left[
\begin{array}{ll}
0_{n \times (m-n)} &  T_{n \times n}  
\end{array}
\right]
\textrm{if $ n \leq m $} 
              & \textrm{or} &
\left[
\begin{array}{cc}
T_{m \times m} \\
0_{(n-m) \times m}
\end{array} 
\right]
\textrm{ if $ n \geq m $} 
\end{array}
\]
where $T_{s \times s}$, for a natural number $s$, is the upper
Toeplitz matrix of order $s$. 
\end{defn}

\begin{lem}
%\marginpar{lem:centralizer}
\label{lem:centralizer}
 Let $n_1, n_2,\ldots, n_l$ be a sequence of natural numbers, such that $n =
  n_1 + n_2 + \cdots + n_l$. Let $A = \oplus_{i=1}^{l} N_{n_{i}}$. 
Then the centralizer, $Z_{\mathrm M_{n}(R)}(A)$ of $A$
  in $\mathrm M_{n}(R)$ consists of matrices of the form 
\[
\begin{bmatrix}
T_{n_{1} \times n_{1}} & T_{n_{1} \times n_2} & \cdots & T_{n_{1}
  \times n_{l}} \\
T_{n_{2} \times n_{1}} & T_{n_{2} \times n_{2}} & \cdots & T_{n_{2}
  \times n_{l}} \\
\vdots & \vdots &   &  \vdots \\
T_{n_{l} \times n_{1}} & T_{n_{l} \times n_{2}} & \cdots & T_{n_{l}
  \times n_{l}} 

\end{bmatrix}
\]
where $T_{n_{i} \times n_{j}}$ for all $i, j$  are rectangular
upper Toeplitz matrices. 
\end{lem}

\begin{defn}(Block Upper Toeplitz Matrix) Let $T_{n_{i} \times n_{j}}$ be a
  rectangular upper Toeplitz matrix of order $n_i \times n_j$, over the
  ring $R$. A matrix of the form
\[
\begin{bmatrix}
T_{n_{1} \times n_{1}} & T_{n_{1} \times n_2} & \cdots & T_{n_{1}
  \times n_{l}} \\
T_{n_{2} \times n_{1}} & T_{n_{2} \times n_{2}} & \cdots & T_{n_{2}
  \times n_{l}} \\
\vdots & \vdots &   &  \vdots \\
T_{n_{l} \times n_{1}} & T_{n_{l} \times n_{2}} & \cdots & T_{n_{l}
  \times n_{l}} 

\end{bmatrix}
\]
is called block upper Toeplitz matrix of order $(n_{1},
  n_{2}, \ldots, n_{l})$ over the ring $R$. 
\end{defn}

In the following lemma we relate the groups of automorphisms
$\G_{\lambda, \mathbf F_q((t))}$ with
the centralizers in $\GL_n(\mathbf F_q)$. 
\begin{lem}
%\marginpar{lem:automorphisms and centralizers}
\label{lem:automorphisms and centralizers}  Let $A = \oplus_{i=1}^{k}N_{\lambda_i}$ and $\lambda = (\lambda_1,
  \lambda_2, \ldots, \lambda_k)$ be a partition. Then the following groups are
  isomorphic, 
\begin{enumerate}
\item The group of automorphisms, $\G_{\lambda, \mathbf F_q((t))} = \rm{Aut}_{\cO}(\cO_{\lambda_1}
  \oplus \cO_{\lambda_2} \oplus \cdots \oplus \cO_{\lambda_k})$. 
\item The centralizer $Z_{\GL_n(\mathbf F_q)}(A)$.
\item The set of invertible block upper Toeplitz matrices of order
  $(\lambda_1, \lambda_2,\ldots, \lambda_k)$ over $\mathbf F_q$. 
\end{enumerate} 
\end{lem}
\begin{pf} Lemma~\ref{lem:centralizer} implies that groups $(2)$ and $(3)$ are
  actually equal. We prove isomorphism between $(1)$ and $(3)$.
  Every $f \in \G_{\lambda}$ can be thought as an invertible matrix of the form 
\[
\left[
\begin{matrix}
f_{11}  &  f_{12}  & \cdots & f_{1k} \\
f_{21}  &  f_{22}  & \cdots & f_{2k} \\
\vdots  &  \vdots  &        &  \vdots \\
f_{k1}  &  f_{k2}  & \cdots &  f_{kk} 
\end{matrix}
\right]
\]
with each $f_{ij} \in \rm{End}_{\cO}(\cO_{\lambda_i}, \cO_{\lambda_j})$. Hence
  it is sufficient to find an isomorphism between 
  $\rm{End}_{\cO}(\cO_{\lambda_i}, \cO_{\lambda_j})$ and rectangular
  Toeplitz matrices of order $\lambda_i \times \lambda_j$ over
  $\mathbf F_q$ which takes composition to matrix multiplication. We prove it only for
  $\lambda_i = \lambda_1$ and $\lambda_j =  \lambda_2$ with $\lambda_1 \geq \lambda_2$, rest of the parts can
  be proved similarly.

Let $\mathcal{T}_{\lambda_1, \lambda_2}$ be the set of rectangular upper
  Toeplitz matrices of order $\lambda_1 \times \lambda_2$ over the field
  $\mathbf F_q$.   
Define a map $\rm{End}_{\cO}(\cO_{\lambda_1}, \cO_{\lambda_2})
  \rightarrow \mathcal{T}_{\lambda_1, \lambda_2}$ by
\[
f \mapsto \left[ \begin{array}{c|c} 0_{(\lambda_1-\lambda_2) \times \lambda_2} & A 
\end{array}
\right]_{\lambda_1 \times \lambda_2}
\]
where
\[
A = \begin{bmatrix}
a_{1} & a_{2} & \cdots & a_{\lambda_2-1} & a_{\lambda_2} \\
0    & a_{1}  & a_{2}&    &    a_{\lambda_2-1} \\
\vdots & \ddots   & \ddots   &  \ddots   &    \vdots \\
      &  \cdots   &       & a_{1}    &  a_{2} \\
0      &   \cdots   &  0    &   0      &   a_{1}   
\end{bmatrix}_{\lambda_2 \times \lambda_2}, 
\]
and elements $a_1, a_2,\ldots, a_{\lambda_2}$ are being determined by the expression $f(1) = a_1 + a_2 \pi
+ \cdots + a_{\lambda_2} \pi^{\lambda_2-1}$. It is straightforward to see that
this map gives the
required isomorphism. 
\end{pf}

%%%%%%%%%%%%%%%%%%%%%%%%%%%%%%%%%%%%%%%%%%%%%%%%%%%%%%%%%%%%%%%%%%%
%
%  Proof of the proposition
%
%
%
%%%%%%%%%%%%%%%%%%%%%%%%%%%%%%%%%%%%%%%%%%%%%%%%%%%%%%%%%%%%%%%%%%
\section{Proof of the Main Theorem}
%\marginpar{section:proof of proposition}
\label{section:proof of proposition}
\noindent  In this section, we present proofs of Proposition~\ref{prop:main} and
Theorem~\ref{thm:main}.\\

Recall that for any character $\rho \in
\hat{K}$, $T(\rho) = \{ g \in \glO{n}{2} \,\,| \,\, \rho^g = \rho
\}$ is the stabilizer of $\rho$ in $\glO{n}{2}$.  By (\ref{eq:orbit}), for each $\psi_A \in \hat{K}$, 
\begin{equation}
\label{lem:inversestab} 
T(\psi_{A}) = \fp^{-1}(Z_{\GL_{n}(\mathcal{O}_1)}(A))
\end{equation}
Fix a section $
s : \mathcal{O}_1 \rightarrow \mathcal{O}_2
$ of the natural quotient map  $\cO_2 \rightarrow \cO_1$,
such that $s (0) = 0$ and $s(1) = 1$. By extending
$s$ entry-wise, we obtain a map $ \s : \mathrm M_{n}(\mathcal{O}_1) \rightarrow \mathrm M_{n}(\mathcal{O}_2).$ 
Observe that restriction of $\s$ to $\glO{n}{1}$ defines a section of $\fp$.
For every matrix $A$ in $\mathrm M_{n}(\mathcal{O}_1)$,
$Z_{\GL_{n}(\mathcal{O}_2)}(\s(A))$ is the centralizer of $\s(A)$ in
$\glO{n}{2}$.
\begin{lem}
%\marginpar{lem:structure}
\label{lem:structure} Assume that $A$ is a split matrix and is in its Jordan canonical form, then
\[
Z_{\GL_{n}(\mathcal{O}_1)}(A) = \fp(Z_{\GL_{n}(\mathcal{O}_2)}(\s(A)))
\]
\end{lem}
\begin{pf}
Let $\alpha = \fp(t)$ for some $t
\in Z_{\glO{n}{2}}(\s(A))$. Then by definition, $t$ satisfies $t \s(A) = \s(A)
t$, which along with the fact
that $\fp$ is a homomorphism, implies  $\fp(t) A = A
\fp(t)$. Hence $\alpha = \fp(t) \in Z_{\GL_{n}(\mathcal{O}_1)}(A)$. This
proves $\fp(Z_{\GL_{n}(\mathcal{O}_2)}(\s(A))) \subseteq
Z_{\GL_{n}(\mathcal{O}_1)}(A)$. 
For the reverse inclusion, since $A$ is a split matrix, by
Theorems ~\ref{thm:primary} and ~\ref{thm:jordan},
\[
A = \oplus_{i=1}^{l} A_{i}
\]
where each $A_{i}$ is a split primary matrix and is of the form 
\[
A_{i} = \oplus_{j=1}^{l_i} J_{\lambda_{ij}}(a_{i}) 
\]
with $a_{i}$'s being distinct elements of the field
$\mathcal{O}_1$. By using $s(0) = 0$ and $s(1) = 1$, we obtain $\s(A_{i})$, for all $i$ is an upper 
triangular matrix with all diagonal
entries equal to $s(a_{i})$ and $\s(A) =
\oplus_{i=1}^{l}(\s(A_i))$. Further for all $i \neq j$, $a_{i} \neq
a_j$ imply that $s(a_{i})-s(a_{j})$ are invertible elements of
the ring $\mathcal{O}_2$. Therefore by Lemma~\ref{lem:generalcentralizers},
%\marginpar{eq:direct}  
\begin{equation}
\label{eq:direct}
Z_{\GL_{n}(\mathcal{O}_2)}(\s(A))  =  \oplus_{i=1}^{l} Z_{\GL_{n_i}(\mathcal{O}_2)}(\s(A_{i}))
\end{equation}
Since we also have 
\[
Z_{\glO{n}{1}}(A) = \oplus_{i=1}^{l} Z_{\GL_{n_{i}(\OO{1})}}(A_{i}),
\]
it is sufficient to prove $Z_{\GL_{n}(\mathcal{O}_1)}(A) \subseteq
\fp(Z_{\GL_{n}(\mathcal{O}_2)}(\s(A)))$ when $A$ is a split primary.

\noindent {\bf Split Primary Case:} 
Now we assume that $A$ is split primary and is in its Jordan canonical form. 
Theorems~\ref{thm:primary} and \ref{thm:jordan} give that $A = aI_{n} + (\oplus_{i=1}^{t}
N_{n_{i}})$ for some $a \in \cO_1$. 
Let $\alpha \in Z_{\GL_{n}(\mathcal{O}_1)}(A)$. By
Lemma~\ref{lem:automorphisms and centralizers}, $\alpha$ is an invertible block Toeplitz matrix of order $(n_1,
n_{2}, \ldots ,n_{t})$ over the ring $\mathcal{O}_1$. Our choice of
section ensures that, 
$\s(A) = s(a)I_n + (\oplus_{i=1}^{t} N_{n_{i}})$, and 
$\s(\alpha)$ is an invertible block Toeplitz matrix of order
$(n_{1},n_{2},\ldots,n_{t})$ over the ring $\mathcal{O}_2$. But then by 
Lemma~\ref{lem:centralizer},
$\s(\alpha) \in Z_{\glO{n}{2}}(\s(A))$. Hence $\alpha = \fp(\s(\alpha)) \in
\fp(Z_{\glO{n}{2}}(\s(A)))$.  
\end{pf}

From the proof of above lemma we obtain, 
\begin{cor}
%\marginpar{cor:structure}
\label{cor:structure}
If $A$ is a split matrix and is in its Jordan canonical form, then $\alpha \in
Z_{\glO{n}{1}}(A)$ if and only if $\s(\alpha) \in Z_{\glO{n}{2}}(\s(A))$.
\end{cor}

%\marginpar{lem:t=kz}   
\begin{cor}
\label{lem:t=kz} If $A$ is a split matrix and is in its Jordan canonical form,
then $ T(\psi_{A}) = K Z_{\glO{n}{2}}(\s(A))$.
\end{cor}
\begin{pf} The inclusion $T(\psi_{A}) \subseteq K Z_{\glO{n}{2}}(\s(A))$
  follows from (\ref{lem:inversestab}) and   
  Lemma~\ref{lem:structure}.      
\end{pf}

\begin{lem}

\label{lem:general diamond lemma}
%\marginpar{lem:general diamond lemma}
Let $G$ be a finite group with two subgroups $N$ and
 $M$, such that $N$ is normal in $G$ and $G = N.M$. If
 $\psi_{1}$ and $\psi_{2}$ are one dimensional representations of $N$ and $M$ respectively such that
 $\psi_{1}(mnm^{-1}) = \psi_{1}(n)$ for all $ m \in M $, $n \in N$ and 
 $\psi_{1}|_{N \cap M} = \psi_{2}|_{N \cap M}$, then
 $\psi_{1}.\psi_{2}$ defined by $\psi_{1}.\psi_{2}(n.m):= \psi_{1}(n)
 \psi_{2}(m)$ is the unique one dimensional representation of
 $G$ extending both $\psi_{1}$ and
 $\psi_{2}$
\end{lem}
\begin{pf} We just prove that $\psi_1.\psi_2$ is well defined, rest of
  the proof is straightforward. Suppose $nm = n'm'$, where $n, n' \in N$ and $m. m' \in M$. Then  $n'^{-1} n = m'm^{-1} \in M \cap N$.
\begin{eqnarray*}
 \psi_{1}.\psi_{2} (nm) & = & \psi_{1}(n) \psi_{2}(m)\\
                       & = & \psi_{1}(n' n'^{-1} n ).\psi_{2}(mm'^{-1}m')
                       \\
 & = & \psi_{1}(n'). \psi_{1}(n'^{-1}nmm'^{-1}).\psi_{2}(m') \\
 & = & \psi_{1}(n').\psi_{2}(m') \\
&  =  & \psi_1\psi_2(n'm') \\
\end{eqnarray*}
\end{pf}

\noindent {\bf Proof of Proposition~\ref{prop:main}}
 It follows from (\ref{eq:orbit})
that orbits of the action of
$\GL_{n}(\mathcal{O}_2)$ on $K$ are the same as orbits of $\mathrm M_{n}(\mathcal{O}_1)$
under the action of $\GL_{n}(\mathcal{O}_1)$, namely the similarity
classes. It is easy to see that, if we can extend the character $\psi_{A}$ from $K$ to $T(\psi_{A})$,
then we can extend any $\psi_{A'}$ in the orbit of $\psi_{A}$ under
the action of $\GL_{n}(\mathcal{O}_2)$ on $T(\psi_{A})$. So to prove the proposition, 
it is enough to choose a representative $A$ of each similarity class of $\mO{n}{1}$ and to extend the
corresponding character $\psi_{A}$ from $K$ to $T(\psi_{A})$. 
We prove existence of this extension in three steps:

\subsection*{Step 1: $A$ is Split Primary} 
Let $A$ be a split primary matrix with unique eigenvalue $a \in \mathcal{O}_1$. Replace $A$ by a matrix in its similarity class of the form  $\oplus_{i=1}^lJ_{\lambda_{i}}(a)$, where each $J_{\lambda_{i}}(a)$ is an elementary Jordan block.

We define a character $\psi_{a} : \mathcal{O}_1 \rightarrow
\mathbb {C}^{*}$ by $ \psi_{a}(x) = \psi(ax)$. The map $ x \mapsto 1 + \pi x$ gives an isomorphism from $\mathcal{O}_1$
onto the subgroup $1 + \pi \mathcal{O}_1$ of the multiplicative group
$\mathcal{O}_2^{*}$. Choose $\chi \in \hat{\cO_2^{*}}$ such that
$\chi(1 + \pi x) = \psi_a(x)$ for all $x\in \cO_1$. Define a character $\tilde{\chi} :
Z_{\GL_{n}(\mathcal{O}_2)}(\s(A)) \rightarrow \mathbb{C}^{*}$ by $\tilde{\chi}(x) = \chi(\det(x))$. 
\begin{lem} The character
  $\tilde{\chi}$ of $Z_{\GL_{n}({\mathcal{O}_2})}(\s(A))$ satisfies 
\[
\tilde{\chi}|_{K \cap Z_{\GL_{n}(\mathcal{O}_2)}(\s(A))} = \psi_{A}|_{K
  \cap Z_{\GL_{n}(\mathcal{O}_2)}(\s(A))}
\] 
\end{lem}

\begin{pf}
By Lemma~\ref{lem:automorphisms and centralizers}, $K \cap Z_{\glO{n}{2}}(\s(A)) = I + \pi
Z_{\mathrm M_{n}(\mathcal{O}_1)}(A)$. If $X = (x_{ij}) \in Z_{\mO{n}{1}}(A)$, then 
by Lemma~\ref{lem:centralizer}, $X$ is a block upper Toeplitz 
matrix. Therefore $\mathrm{Tr}(AX) = a(x_{11}+ x_{22} + \cdots + x_{nn})$. 
We have
\begin{eqnarray*}
\psi_{A}(I + \pi X) & = & \psi (\tr (AX)) \\
                  & = & \psi (a (x_{11} + x_{22}
+ \cdots + x_{nn})) \\
                  & = & \chi (\det(I + \pi X)) = \tilde{\chi}(I + \pi X)
\end{eqnarray*}    
\end{pf}

Applying Lemma~\ref{lem:general diamond lemma} to the group
$T(\psi_A)$ with its subgroups  $K$ and  $Z_{\glO{n}{2}}(\s(A))$, and characters $\psi_1 = \psi_A$ and
$\psi_2 = \tilde{\chi}$ we obtain that the character
$\psi_{A}.\tilde{\chi}$ is an extension of $\psi_A$ from $K$ to $T(\psi_{A})$.  

\subsection*{Step 2:  $A$ is split }
Let $A$ be a split matrix with distinct eigenvalues $a_{1}, a_{2},
\ldots ,a_{l}$. Then by Theorem~\ref{thm:primary}, $A$ can be written as $\oplus_{i=1}^{l} A_{i}$, where
each $A_{i}$ is a split primary matrix, say of order $n_{i}$, and has a unique
eigenvalue $a_{i}$. We may assume each $A_{i}$ in its
Jordan canonical form. Then by (\ref{eq:direct}),
\[
K \cap Z_{\glO{n}{2}}(\s(A)) = \oplus_{i=1}^{l} (K \cap
Z_{\GL_{n_i}(\mathcal{O}_2)}(\s(A_{i}))) 
\]

As in the Step 1, define the characters $\tilde{\chi_{i}}$ of
$Z_{\GL_{n_i}(\mathcal{O}_2)}(\s(A_i))$ such that 
\[
\tilde{\chi_{i}}|_{K \cap
Z_{\GL_{n_i}(\mathcal{O}_2)}(\s(A_{i}))} =
    \psi_{A_{i}}|_{K \cap
Z_{\GL_{n_i}(\mathcal{O}_2)}(\s(A_{i}))}
\]
Then the character $\tilde{\chi} = \tilde{\chi}_{1} \times \tilde{\chi}_{2}
\times \ldots \tilde{\chi}_{l}$ is a character of $Z_{\glO{n}{2}}(\s(A))$, such that
\[
\tilde{\chi}|_{K \cap Z_{\GL_{n}(\mathcal{O}_2)}(\s(A))} = \psi_{A}|_{K
  \cap Z_{\GL_{n}(\mathcal{O}_2)}(\s(A))}
\]
Again by Lemma~\ref{lem:general diamond lemma},  $\psi_{A}.\tilde{\chi}$
is an extension of $\psi_{A}$ from $K$ to $T(\psi_{A})$. 

\subsection*{Step 3 : General case} 
Let $\tilde{\mathcal{O}_1}$ be a splitting field for the characteristic
 polynomial of $A$ and let $\tilde {\mathcal{O}_2}$ be the corresponding unramified extension
 of $\mathcal{O}_2$. Let $\tilde{K}= \rm{Ker}(\GL_n(\tilde{\cO_2}) \rightarrow \GL_n(\tilde{\cO_1}))$ under the natural 
quotient map, and $\tilde{\psi} : \tilde{\cO_1} \rightarrow \mathbb C^{*}$ be a character 
such that $\tilde{\psi}|_{\cO_1} = \psi$. 
Then $\tilde{\psi}_{A} :
\tilde{K} \rightarrow  \mathbb C^{*},$ defined by 
\[
\tilde{\psi}_{A} (I + \pi X) = \tilde{\psi}(\tr (AX))
\]
 is a character of $\tilde{K}$. Let $\tilde{T}(\psi_{A})$ be
 the stabilizer of $\tilde{\psi}_{A}$ in $\GL_{n}(\tilde{\cO_2})$. Since $A$
 splits over $\tilde{\cO_1}$, by Step~2, there exists a character
 $\tilde{\chi} : \tilde{T}(\psi_{A}) \rightarrow \mathbb C^{*}$, such that
\[
\tilde{\chi}|_{\tilde{K}} = \tilde{\psi}_{A}.
\] 
Define a character $\chi:
 T(\psi_{A}) \rightarrow  \mathbb C^{*}$ by $\chi =
 \tilde{\chi}|_{T(\psi_A)}$. Then $\chi$ is an extension of
 $\psi_{A}$ to $T(\psi_{A})$. This completes the proof of
 Proposition~\ref{prop:main}. \\

Fix an extension $\chi_A$ of $\psi_A$ from $K$ to $T(\psi_A)$ and let $\mathcal{S}$ denote 
the set of similarity classes
of $\mO{n}{1}$. By $(\ref{lem:inversestab})$, the groups $T(\psi_A)/K$ and $Z_{\glO{n}{1}}(A)$ are isomorphic. Therefore by Clifford Theory, there exists a bijection between the sets 
\begin{eqnarray}
\label{correspondence}
 \coprod_{A \in \mathcal{S}} \{ \mathrm{Irr}(Z_{\glO{n}{1}}(A))\} &
\longleftrightarrow & \rm{Irr}(\GL_{n}(\cO_2)),
\end{eqnarray}
given by, 
\begin{eqnarray}
\phi \mapsto \mathrm{Ind}_{T(\psi_A)}^{\glO{n}{2}}(\chi_{A} \otimes \phi).
\end{eqnarray}

As $[\glO{n}{2} : T(\psi_A)]$
= $[\glO{n}{1} : Z_{\glO{n}{1}}(A)]$, this already proves that there exists a dimension preserving bijection between the
sets $\rm{Irr}(\glO{n}{2})$ and $\rm{Irr}(\GL_n(\cO'_2))$.
 To prove that this
bijection is canonical we need to do little more work.

Let $\cO'$ be an another ring of integers of local non-Archimedean field $F'$, such that residue fields of both
$\cO$ and $\cO'$ are isomorphic. We fix an isomorphism $\phi$
between their residue fields. From now onwards we shall assume that section $
s : \mathcal{O}_1 \rightarrow \mathcal{O}_2
$ satisfies $s (0) = 0$ and $s|_{\cO_1^{*}}$ is multiplicative. The existence
and uniqueness of this section is proved in, for example, Serre~\cite[Prop. 8]{MR0354618}. In the sequel, this unique section will be called \emph{the multiplicative section} of $\cO_1$ (or of $\cO^{*}_1$)
(depending on the domain). Given above isomorphism $\phi$,  
\begin{lem} 
\label{lem:canonical isomorphism} There exists a canonical isomorphism between groups $\hat{\cO_2^{*}}$
  and $\hat{\cO_2^{'*}}$. 
\end{lem}
\begin{pf} Let $s : \cO^{*}_1 \rightarrow \cO_2^{*}$ and $s': \cO^{'*}_1 \rightarrow
  \cO_2^{'*}$ be the multiplicative sections of $\cO^*_1$ and $\cO^{'*}_1$ respectively. Then the following exact sequences split,
\[
\xymatrix{
0 \ar[r] & \cO_1  \ar[r]_{i} & \cO_2^{*} \ar[r] & \cO^{*}_1 \ar[r] \ar@/^/@{.>}[l]^s & 1 \\
0 \ar[r] & \cO'_1 \ar[r]_{i'} & \cO_2^{'*} \ar[r] & \cO^{'*}_1 \ar[r] \ar@/_/@{.>}[l]_{s'} & 1 }
\]  
The uniqueness of the sections $s$ and $s'$ implies the existence of a unique isomorphism $f :\cO_2^{*} \rightarrow \cO_2^{'*}$
such that
$f \circ  s = s' \circ \phi$ and $f \circ i = i' \circ \phi$. This
gives a canonical isomorphism between
$\hat{\cO_2^{*}}$ and $\hat{\cO_2^{'*}}$.
\end{pf} 

\noindent Let $K' = \rm{Ker} (\GL_n(\cO'_2) \rightarrow \GL_n(\cO'_1))$. Then 
\[
K \cong K' \cong \mO{n}{1}
\] 
hence the set $\{\psi_A\,\,| \,\, A \in \mO{n}{1}\}$ can also be
thought as the set of
characters of $K'$. Let $T'(\psi_A)$ denote the stabilizer of character $\psi_A$
in $\GL_n(\cO'_2)$. By (\ref{lem:inversestab}), groups $T'(\psi_A)/K'$ and $T(\psi_A)/K$ are
canonically isomorphic. Further to prove that there exists a canonical
bijection between $\rm{Irr}(\glO{n}{2})$ and
$\rm{Irr}(\GL_n(\cO'_2))$, it is sufficient to prove that for a given 
$A \in \mO{n}{1}$, and an extension 
$\chi_{A}: T(\psi_A) \rightarrow \mathbf C^{*}$ of $\psi_A$, there exists a canonical extension
$\chi'_A: T'(\psi_A) \rightarrow \mathbf C^{*}$ of $\psi_A$ from $K'$
to $T'(\psi_A)$. For that, in steps 1 and 2 of the proof of Proposition~\ref{prop:main}, choose the
character of $\cO_2^{'*}$ by using the given character of
$\cO_2^{*}$ and canonical isomorphism between $\hat{\cO_2^*}$
and $\hat{\cO_2^{'*}}$; the rest of the argument follows easily for these steps. 
For step 3, any isomorphism between
$\tilde\cO_1$ and $\tilde\cO'_1$ that extends $\phi:\cO_1\to\cO'_1$ is defined uniquely up to an element of the Galois
group. To complete the proof for this step, we 
observe that, if $\tilde{\psi}$ is an extension of $\psi$ from $\cO_1$ to
$\tilde{\cO_1}$ and $\gamma$ is an element of Galois group
$\rm{Gal}(\tilde{\cO_1}/\cO_1)$ then by
definition $(\tilde{\psi} \circ \gamma)_{A} = \tilde{\psi}_A \circ
\gamma$ where on the right side element $\gamma$ is thought as scalar matrix
with all its diagonal entries equal to $\gamma$. Choose
$\tilde{\chi}_A$ an extension of $\tilde{\psi_A}$ from $\tilde{K}$ to
$\tilde{T}(\psi_A)$ and let $\tilde{\gamma}$ be a lift of $\gamma$
from $\rm{Gal}(\tilde{\cO_1}/\cO_1)$ to $\rm{Gal}(\tilde{F}/F)$, which
takes the maximal ideal of $\tilde\cO$ to itself (for existence of this see \cite[p. 26]{MR0215665}), then $\tilde{\chi} \circ
\tilde{\gamma}$ (again $\tilde{\gamma}$ is
thought as scalar matrix with diagonal entries equal to
$\tilde{\gamma}$) extends $ \tilde{\psi}_A \circ
\gamma$. The restrictions of $\tilde{\chi}_A$ and $\tilde{\chi}_A
\circ \tilde{\gamma}$ to $T(\psi_A)$ coincide. This completes the proof of Theorem~\ref{thm:main}.
 
\begin{cor} The isomorphism type of group algebra $ \mathbf C[\glO{n}{2}]$ depends 
only on the cardinality of the residue field of $\cO$.
\end{cor}

This can be restated as,
\begin{cor}
The number and dimensions of irreducible representations of groups
$\glO{n}{2}$ depends only on the residue field. 
\end{cor}
The equivalence between number of conjugacy classes and irreducible
representations further gives,
\begin{cor} The number of conjugacy classes of groups $\glO{n}{2}$ depend
  only on the cardinality of residue fields.
\end{cor}
\begin{remark}
  In a forthcoming paper, we will sharpen this result by showing that the class equation of $\glO n2$ depends only on the cardinality of the residue field.
\end{remark}

\section{Complexity of the problem}
\label{complexity}
In this section we comment on the complexity of the
problem of constructing all the irreducible representations of
$\glO{n}{2}$.\\

In this context, Aubert-Onn-Prasad-Stasinski
have proved the following (\cite[Theorem~6.1]{AOPS})
\begin{thm}
Let $F = \mathbf F_q((t))$ be a local function field. Then the problem of
constructing all the irreducible representations of the following are
equivalent:
\begin{enumerate}
\item $G_{2^n, F}$ for all $n \in \mathbf N$.
\item $G_{k^n, F}$ for all $k, n \in \mathbf N$. 
\item $G_{\lambda, E}$ for all partitions $\lambda$ and all unramified
  extensions $E$ of $F$.
\end{enumerate}   
\end{thm}

\noindent The above, combined with Theorem~\ref{thm:main}, proves the following, 

\begin{thm} Let $\cO$ be ring of integers of a non-Archimedean local field
  $F$, such that the residue field has cardinality $q$. 
  Then the problem of constructing irreducible representations of the
  following groups are equivalent:
\begin{enumerate}
\item $\glO{n}{2}$ for all $n \in \mathbf N$.
\item $G_{\lambda, E}$ for all partitions $\lambda$ and all unramified
  extensions $E$ of $\mathbf F_q((t))$.   
\end{enumerate}
\end{thm}

%%%%%%%%%%%%%%%%%%%%%%%%%%%%%%%%%%%%%%%%%%%%%%%%%%%%%%%%
%
% Applications 
%%%%%%%%%%%%%%%%%%%%%%%%%%%%%%%%%%%%%%%%%%%%%%%%%%%%%%%%%%%%%%%5
\section{Applications}
\label{sec:applications}
%\marginpar{sec:applications}

In this section, we discuss a few applications of the theory developed so
far. In particular, we discuss the relation between the representation zeta
function of $\glO{n}{2}$ and those of centralizers in $\glO{n}{1}$.
We also construct all the irreducible representations of groups $\glO{2}{2}$,
$\glO{3}{2}$, $\glO{4}{2}$ and obtain their representation zeta functions. 

Recall the following definition from Section~1, 
\begin{defn}(Representation Zeta function)
Let $G$ be a group. The representation zeta function of $G$ is the function
\[
R_{G}({\D}) = \sum_{\rho \in \irr{G}} \D^{\dim{\rho}} \in \mathbb Z[\mathcal{D}]
\] 
\end{defn}
\subsection{Representation zeta function of $\glO{n}{2}$}
Let $\mathcal{S}$ be the set of similarity classes of $\mO{n}{1}$.  
From (\ref{correspondence}) it is clear that 
representations of centralizers play an important role in determining
irreducible representations of $\glO{n}{2}$. Moreover we obtain the
following relation between their representation zeta functions.
\begin{equation}
\label{eq:similarity_zeta}
R_{\glO{n}{2}}(\D) = \sum_{A \in \mathcal{S}}
R_{Z_{\mO{n}{1}}(A)}(\D^{[\glO{n}{1} : Z_{\glO{n}{1}}(A)]})
\end{equation}
where, $[\glO{n}{1} : Z_{\glO{n}{1}}(A)] = [ \glO{n}{2} : T(\psi_A)]$ is the index of
  $Z_{\glO{n}{1}}(A)$ in $\glO{n}{1}$. Following Green~\cite{MR0072878}, a similarity class $c$ 
of $\mathrm M_n(\cO_1)$ can be denoted by the
symbol
\[
c = (\ldots,f^{\nu_c(f)},\ldots)
\]
where $f$ is an irreducible polynomial appearing in the characteristic
polynomial of $c$ and $\nu_c(f)$ is the partition associated with $f$
in the canonical form of $c$. 

Let $c = (\ldots, f^{\nu_c(f)}, \ldots)$. Let $d$ be a positive integer, and
let $\nu$ be a partition other than zero. Let $r_c(d, \nu)$ be the
number of $f$ appearing in the characteristic polynomial of $c$ with
degree $d$ and $\nu_c(f) = \nu$. Let $\rho_c(\nu)$ be the partition 
\[
\{ n^{r_c(n, \nu)}, (n-1)^{r_c(n-1, \nu)},\ldots\}
\]
Then two classes $b$ and $c$ are of the same {\em type} if and only if
$\rho_b(\nu) = \rho_c(\nu)$ for each non zero partition $\nu$. By
abusing notation we shall also say matrices of class $c$ and $d$ have
same type.   

Let $\rho_{\nu}$ be a partition-valued function on the nonzero
partitions $\nu$ ($\rho_{\nu}$ may take value zero). The condition for
$\rho_{\nu}$ to describe a type of $\mathrm M_n(\cO_1)$ is 
\[
\sum_{\nu} |\rho_{\nu}||\nu| = n 
\]
The total number $t(n)$ of functions $\rho_{\nu}$ satisfying above expression
is independent of $q$, and so is the number of types of
$\mathrm M_n(\cO_1)$ (for large enough $q$).
The following lemma (which is easy) underlines the importance of types in the calculation of the representation zeta functions of the groups $\glO n2$:

\begin{lem} If matrices $A$ and $B$ in $\mO{n}{1}$ are of same {\em
type} then their centralizers are isomorphic. 
\end{lem}

Let $\mathcal{T}$ denote the set of representatives of types of $\mO{n}{1}$ and 
for each $A \in \mathcal{T}$, let $n_{A}$ be the total number of similarity classes
of type $A$. The expression (\ref{eq:similarity_zeta}) simplifies to
\begin{equation}
\label{eq:type_zeta}
R_{\glO{n}{2}}(\D) = \sum_{A \in \mathcal{T}}
n_{A} R_{Z_{\mO{n}{1}}(T)}(D^{[\glO{n}{1} : Z_{\glO{n}{1}}(A)]})
\end{equation}
Summarising the discussion so far, to determine the irreducible
representations of groups $\glO{n}{2}$, it is sufficient to determine
the representations of the centralizers $Z_{\glO{n}{1}}(A)$ where $A$
varies over the set of types of $\mathrm M_n(\cO_1)$. But determining
representations of groups $Z_{\glO{n}{1}}(A)$ for general $n$ is still an open
problem. We discuss representations of these groups for $n = 2$, $n= 3$, and $n= 4$.

We shall use the following theorems in the sequel. For proofs of these results see for example
\cite[Chapter~1]{MR0201472}. 
\begin{thm}
%\marginpar{lem:field}
\label{lem:min=characteristic}
%\marginpar{lem:min=characteristic} 
Let $A \in \mathrm{M}_n(\mathbf F_q)$, then
  centralizer of $A$ in $\mathrm M_{n}(\F_{q})$ is the
  algebra $\mathbf F_{q}[A]$ if and only if the minimal polynomial
  of $A$ coincides with its characteristic polynomial. Moreover, in this case,
  $\mathrm{dim}_{\mathbf F_{q}}Z_{\mathrm M_{n}(\F_{q})}(A) = n$.  
\end{thm}
\begin{thm}
\label{thm:field}
%\marginpar{thm:field} 
Let $A \in \mathrm M_n(\mathbf F_q)$, then its centralizer in $\mathrm M_n(\mathbf
F_q)$ is a field if and only if its characteristic polynomial is
irreducible over $\mathbf F_{q}$.  
\end{thm}
\subsection{Representation zeta function  of $\glO{2}{2}$}
\label{subsection:representations of 22}
%\marginpar{subsection:representations of 22}
The irreducible
representations of groups $\glO{2}{2}$ are already described by
Nagornyi \cite{MR0491919} and Onn \cite {MR2456275}. Since it falls out of
our discussion very easily and is used in representation theory of
groups $\glO{4}{2}$, we add its brief description also. For representation theory of
groups $\glO{n}{1}$, we refer to Green~\cite{MR0072878} and
Steinberg~\cite{MR0041851}. In Table\ref{Table1} we describe types of
$\mathrm M_{2}(\cO_1)$ (set of $2 \times 2$ matrices over $\cO_1$) with their centralizers.  To determine centralizers, wherever required, we have used 
Theorems~\ref{lem:min=characteristic} and ~\ref{thm:field}.  

\begin{table}
The element $\rho$ and $\sigma$ are primitive elements of $\F_q$ and $\F_{q^2}$
respectively, such that $\rho = \sigma^{q+1}$.
\[
\scriptsize{
\begin{array}{|c|c|c|c|}
\hline
 Type & \mathrm{ Number\; of \; similarity } &
     \mathrm{Isomorphism\; type \; } & Indices \\
   &   \mathrm{classes \; of \;given}   &  \mathrm{of \; centralizer}
   &  \\
A   &  \mathrm{type} \;\; (n_{A}) & Z_{\glO{2}{1}}(A)   & \\
\hline

\left( \begin{array}{cc}
 \rho^a& 0 \\
0 & \rho^a 
\end{array} \right)      &   q   &   R_{\glO{2}{1}}(\D) & 1 \\
&&&\\
\left( \begin{array}{cc} 
\rho^a & 0 \\
0 & \rho^b \\
\end{array}
\right)      &  \frac{1}{2}q(q-1) & (q-1)^2\D & q(q+1) \\
&&&\\      
\left(
\begin{array}{cc}
\rho^a & 1 \\
0 & \rho^a \\
\end{array}
\right)                  &  q  & q(q-1) \D & q^2-1 \\
&&&\\
\left(
\begin{array}{cc}
\sigma^a &  0 \\
0   &  \sigma^{aq} 
\end{array}
\right)                 &   \frac{1}{2}q(q-1) & (q^2 -1) \D & q^2-q \\
\hline
\end{array}
}
\]
\caption{\label{Table1} Group $\glO{2}{2}$}
\end{table}
The representation zeta function of the group
$\glO{2}{1}$ (see
Steinberg~\cite{MR0041851}) is
\begin{eqnarray*}
R_{\glO{2}{1}}(\D) & = & (q-1) \D + (q-1) \D^{q} + \frac{1}{2}(q-1)(q-2)
\D^{q+1} \nonumber \\
&   &  + \frac{1}{2} q(q-1) \D^{q-1}.
\end{eqnarray*} Feeding all this data into
\ref{eq:type_zeta}, we easily obtain the
representation zeta function of $\glO{2}{2}$:
\begin{eqnarray}
R_{\glO{2}{2}}(\D) &  = &  q R_{\glO{2}{1}}(\D) + \frac{1}{2} q(q-1)^3
\D^{q(q+1)} + q^2 (q-1) \D^{q^2-1} \nonumber \\ 
  &   &  + \frac{1}{2}q(q+1)(q-1)^2 \D^{q^2-q} 
\end{eqnarray}

\subsection{Representation zeta function of $\glO{3}{2}$} 
Partial results regarding representations of groups $\glO{3}{2}$ are
already given by Nagornyi~\cite{MR0498881}. We complete his results
for these groups. 

In Table~\ref{Table2} we describe types and their corresponding
centralizers for the group $\mO{3}{1}$. 
\begin{table}
 The elements $\rho$, $\sigma$ and $\tau$ are primitive elements of $\F_q$,
$\F_{q^2}$ and $\F_{q^3}$ respectively, such that $\rho = \sigma^{q+1} = \tau^{q^2
  + q +1}$ and $\sigma = \tau^{q^2 + 1}$.
\[
\scriptsize{
\begin{array}{|c|c|c|c|}
\hline
 Type & \mathrm{ Number\; of \; similarity } &
     \mathrm{Isomorphism\; type \; } & Index \\
   &   \mathrm{classes \; of \;given}   &  \mathrm{of \; centralizer}
   & [Z_{\glO{3}{1}}(A) : \glO{3}{1}] \\
A   &  \mathrm{type} \;\; (n_{A}) & Z_{\glO{3}{1}}(A)   & \\
\hline
\left(
\begin{array}{ccc}
\rho^a & 0 & 0 \\
0 & \rho^a & 0 \\
0 & 0 & \rho^a 
\end{array}
\right)             & q  &  \mathcal{R}_{\glO{3}{1}}(\D) & 1 \\
&&& \\
\left( 
\begin{array}{ccc}
\rho^a & 0 & 0 \\
0 & \rho^a & 0 \\
0 & 0 & \rho^b 
\end{array}
\right)            &  q(q-1)  &  (q-1) \mathcal{R}_{\glO{2}{1}}(\D) &
q^2(q^2+q+1) \\
&&&\\

\left( 
\begin{array}{ccc}
\rho^a & 0 & 0 \\
0 & \rho^b & 0 \\
0 & 0 & \rho^c 
\end{array}
\right)           &  \frac{1}{6}q(q-1)(q-2)   &  (q-1)^3 \D &
q^3(q+1)(q^2+q+1) \\
&&&\\

\left(
\begin{array}{ccc}
\rho^a & 1 & 0 \\
0 & \rho^a & 0 \\
0 & 0 & \rho^b 
\end{array}
\right)           &  q(q-1)  &  q (q-1)^2 \D & q^2(q^3-1)(q+1) \\
&&&\\
\left(
\begin{array}{ccc}
\rho^a & 1 & 0 \\
0 & \rho^a & 0 \\
0 & 0 & \rho^a 
\end{array}
\right)          &  q   & \mathcal{R}_{G_{(2,1)}}(\D) & (q^3-1)(q+1)
\\
&&&\\ 
\left(
\begin{array}{ccc}
\rho^a & 1 & 0 \\
0 & \rho^a & 1 \\
0 & 0 & \rho^a 
\end{array}
\right)         &  q  &   q^2(q-1) \D & q(q^3-1)(q^2-1) \\
&&&\\
\left(
\begin{array}{ccc}
\rho^a & 0 & o \\
0 & \sigma^{a} & 0 \\
0 &   0 & \sigma^{aq}\\
\end{array}
\right)          &  \frac{1}{2}q^{2}(q-1)  &  (q-1)(q^2-1) \D &
q^3(q^3-1) \\
&&&\\
\left(
\begin{array}{ccc}
\tau & 0 & 0 \\
0 & \tau^{bq} & 0 \\
0 & 0 & \tau^{bq^2} 
\end{array}
\right)           &          \frac{1}{3} q(q^2 -1)  &   (q^3 - 1) \D &
q^3(q-1)^2(q+1) \\
\hline
\end{array}
}
\]
\caption{\label{Table2} Group $\glO{3}{2}$}
\end{table}
The representation zeta
function of $\glO{3}{1}$ (Steinberg~\cite{MR0041851}) is 
\begin{eqnarray}
\mathcal{R}_{\glO{3}{1}}(\mathcal{D}) & = & (q-1) \D + (q-1) \D^{q^2 + q} + (q-1)
\D^{q^3}  \nonumber \\
   &  &  (q-1)(q-2) \D^{q^2 + q + 1} + (q-1)(q-2) \D^{q(q^2 + q + 1)}
  \nonumber \\
 &  &  + \frac{1}{6}(q-1)(q-2)(q-3) \D^{(q+1)(q^2 + q + 1)} \nonumber \\ 
&  &  + \frac{1}{2} q(q-1)^2 \D^{(q-1)(q^2 + q + 1)} \nonumber \\
&  & + \frac{1}{3} q(q-1)(q+1) \D^{(q+1)(q-1)^{2}} 
\end{eqnarray} 
The irreducible representations of
all the centralizers  appearing in Table~\ref{Table2} except $G_{2,1}$ are
either very easy or
well known. Onn~\cite[Theorem 4.1]{MR2456275} has described all the 
irreducible representations of groups $G_{(\ell,1)}$ for $\ell > 1$.
As a consequence, we have

\begin{lem} 
\label{lem:G(2,1)} 
%\marginpar{lem:G(2,1)}
The representation zeta function of the group $G_{(2,1)}$ is 
\[
R_{G_{(2,1)}}(\mathcal{D}) =  (q-1)^2 \mathcal{D} + (q^2 - 1)\mathcal{D}^{q-1}
+ (q - 1)^{3} \mathcal{D}^{q}
\]
\end{lem}
Collecting all the pieces together, we obtain the expression for representation zeta function of $\glO{3}{2}$:
\begin{eqnarray}
\mathcal{R}_{\glO{3}{2}}(\mathcal{D}) & =  & q
\mathcal{R}_{\glO{3}{1}}(\mathcal{D}) + q(q-1)^2
\mathcal{R}_{\glO{2}{1}}(\mathcal{D}^{q^2(q^2+q+1)}) \nonumber \\ 
&  &  + \frac{1}{6}q(q-2)(q-1)^4
\mathcal{D}^{q^3(q+1)(q^2+q+1)} + \nonumber \\ 
&  &  q^2(q-1)^3 \mathcal{D}^{q^2(q^3-1)(q+1)} + q
\mathcal{R}_{\G_{(2,1)}}(\mathcal{D}^{(q^3-1)(q+1)}) \nonumber \\
&  &  + q^3 (q-1) \mathcal{D}^{q(q^3-1)(q^2-1)} \nonumber 
\end{eqnarray}
\begin{eqnarray}
& & + \frac{1}{2}q^2(q-1)^2(q^2-1)\mathcal{D}^{q^3(q^3-1)} \nonumber \\
&  &  + \frac{1}{3}q(q^2-1)(q^3-1)
\mathcal{D}^{q^3(q-1)^2(q+1)} 
\end{eqnarray}

\subsection{Representation zeta function of $\glO{4}{2}$}
In this section we discuss
representation theory of groups $\glO{4}{2}$.

In Table~\ref{Table3}, we give all the data required for the
representations of $\glO{4}{2}$. The expression for representation
zeta function of $\glO{4}{1}$ is rather long, so we omit the details
here (see Steinberg~\cite{MR0041851}). Among the other
centralizers appearing in this Table only the results regarding the
representations of group $G_{(2,1,1)}$ are not clear from our
discussion so far. We follow a method of Uri Onn to discuss representations of
these groups.  

For the proof of next Proposition (which follows from the theory of finite Heisenberg groups) we refer
Bushnell-Fr{\"o}hlich~\cite[Prop 8.3.3]{MR701540}.
\begin{prop}
%\marginpar{prop:Heisenberg like groups}
\label{prop:Heisenberg like groups} Let $1 \rightarrow N \rightarrow G \rightarrow^{\phi} G/N \rightarrow 1$ be an extension, where $V = G/N$ is an elementary finite abelian $p$-group
  so also viewed as a finite dimensional vector
  space over $\F_p$. Further, let $\chi: N \mapsto \mathbb C^{*}$ be a nontrivial
  character such that $G$ stabilizes $\chi$. Assume furthermore that
  $h_{\chi}(g_1 N, g_2 N) =  \langle g_1 N, g_2 N \rangle_{\chi} = \chi([g_1 , g_2 ])$ is an alternating
  nondegenerate bilinear form on $V$. Then there exists a unique irreducible
  representation $\rho_{\chi}$ of $G$ such that $\rho_{\chi}|_{N}$ is
  $\chi$-isotypic. Moreover, $\rm{dim}(\rho_{\chi})^2 = [G : N] $. 
\end{prop} 

\begin{lem}
\label{lem:G(2,1,1)}
The representation zeta function of the group $G_{(2,1,1)}$ is 
\[
\mathcal{R}_{G_{(2,1,1)}}(\mathcal{D})  =  (q-1)^2
\mathcal{R}_{\glO{2}{1}}(\mathcal{D}^{q^2}) + (q-1) \mathcal{R}_{(\cO_1^{2} \times
  \cO_{1}^{2}) \rtimes \G_{(1,1)}}(\mathcal{D}),
\]
where 
\begin{eqnarray}
\mathcal{R}_{(\cO_1^{2} \times
  \cO_1^{2}) \rtimes \G_{(1,1)}}(\mathcal{D}) & =  & 
\mathcal{R}_{(1,1)}(\mathcal{D}) + 2(q-1) \mathcal{D}^{q^2 -1} + (q-1)^2
  \mathcal{D}^(q^2-1)q \nonumber \\
&  &  + (q+2) \mathcal{D}^{(q^2-1)(q-1)}.
\end{eqnarray}
\end{lem}

\begin{pf}
Using the notation in the proof of Lemma~\ref{lem:automorphisms and centralizers},
let $H$ be the kernel of map $\G_{(2,1,1)} \rightarrow \G_{1} \times
\G_{(1,1)}$,
\[
g \mapsto \left(g_{11} (\rm{mod}\wp), \left( \begin{matrix} g_{22} & g_{23} \\ g_{32} & g_{33}
\end{matrix} \right)\right).
\]
Then
\[
H = I + \left[
\begin{matrix}
\wp & \wp & \wp \\
\cO_1 &   &    \\
\cO_1  &   &  
\end{matrix}
\right] 
\]
The centre of H, i.e. $Z(H) \cong \cO_1$. Firstly we claim that $H$ has $q-1$ irreducible representations of dimension $q^2$ which lie above the non-trivial characters of $Z(H)$. We identify $Z(H)$ with its dual by $z \mapsto \psi_{z}(.) =
\psi(\tr(z.))$. $H$ stabilizes these characters of $Z(H)$ and
furthermore, each of the non-trivial characters gives rise to an alternating
non-degenerate bilinear form $\langle h_1 Z(H), h_2 Z(H)
\rangle_{\psi_{z}}= $ $\psi(\tr(z[h_1,
  h_2]))$ on $H/Z(H)$. Proposition~\ref{prop:Heisenberg like groups} gives $q-1$ pairwise inequivalent irreducible representations of dimension $|H/Z(H)|^{1/2} =
q^2$. This proves the claim. Furthermore, the group $\G_{2,1,1}$ stabilizes each of these
representations of $H$. Let $\rho_{\chi} \in \widehat{H}$ be such a representation
lying over a non-trivial character $\chi \in \widehat{Z(H)}$. We claim that the
representation $\rho_{\chi}$ can be extended to $\G_{(2,1,1)}$.
Let $H^i = Z(H) \times \cO_1 \times \cO_1 $ be the pre-image in $H$ of the maximal isotropic subgroup $\cO_1\times \cO_1$ for the above bilinear form. 
Let $\chi^{i}$ be any extension of $\chi$ to $H^i$.
Indeed the subgroup $G_{2} \times
G_{(1,1)}$ stabilizes both $H^{i}$ and $\chi^i$. Let $\tilde{\chi}$ be an
extension of $\chi$ to $G_{2} \times G_{(1,1)}$. Then by
Lemma~\ref{lem:general diamond lemma}, $\chi^i.\tilde{\chi}
(a.b) = \chi^{i}(a) \tilde{\chi}(b)$ for all $a \in H^i$ and $b \in G_{2}
\times G_{(1,1)}$ is a well defined linear character of $H^i.(G_{2} \times
G_{(1,1)})$. By the proof of Proposition~\ref{prop:Heisenberg like groups}, $\rho_{\chi}$ does not 
depend on the choice of isotropy group and the extension $\chi^i$. Therefore for the induced representation 
\[
\rho_{\chi^i} = \rm{ind}_{H^i(G_{2} \times
G_{(1,1)})}^{H.G_{2} \times G_{(1,1)}} (\chi^i),
\]
$\rho_{\chi} \leq \rho_{\chi^i}$, and as $\rm{dim} \rho_{\chi} =
\rm{dim} \rho_{\chi^i} = q^2 $ we conclude that $\rho_{\chi^i}$ is an
extension of $\rho_{\chi}$ to $\G_{(2,1,1)}$. By the Clifford theory, it follows that all the
representations of $\G_{(2,1,1)}$ which lie above $\rho_{\chi}$ are of the
form $\{ \rho_{\chi^i}. \phi \, | \, \phi \in \G_{(2,1,1)}/H \}$. Hence the
contribution to representation zeta function of $\G_{(2,1,1)}$ from these
representations is $(q-1) \mathcal{R}_{\cO_{1}^{*} \times \GL_2(\cO_1)}
(\mathcal{D}^{q^2})$. 
The remaining representations correspond to representations of $H$ whose
central character is trivial, that is, representations pulled back from
$((\cO_1^{2} \times \cO_1^2) \rtimes  G_{(1,1)}) \times \cO_{1}^{*}$. The action of
$G_{(1,1)}$ on $\cO_1^2 \times \cO_1^2$ is given by
\[
\left(
\begin{matrix}
1 &  \\
 &  D 
\end{matrix}
\right)
\left(
\begin{matrix}
1 & \pi v \\
w & I  
\end{matrix}
\right)
\left(
\begin{matrix}
1 &  \\
  & D^{-1} 
\end{matrix}
\right)
= 
\left(
\begin{matrix}
1  &  \pi v D^{-1} \\
Dw  &  I  
\end{matrix}
\right),
\]
After a choice of identification of $\cO_1^2 \times \cO_1^2$ with its dual: $\langle (\hat{v}, \hat{w}),
(v, w) \rangle = \psi (v \hat{v} + w \hat{w})$, we get
\[
g^{-1}(\hat{v}, \hat{w}) = (D^{-1} \hat{v}, \hat{w} D),\; \mathrm{where} \;\; g = \left(
\begin{matrix}
1 &  \\
 &  D 
\end{matrix}
\right),  
\]
and the orbits and stablizers of this action are given by 
are given by 
\begin{displaymath}
\begin{array}{|c|c|c|}
\hline
   &  \rm{Orbits}   & \rm{Stabilizers} \\
\hline
(1) & \left[\begin{matrix} \left(\begin{matrix} 0 \\ 0 \end{matrix} \right), &
      \left(\begin{matrix} 0 & 0 \end{matrix} \right) \end{matrix} \right] &
  \G_{(1,1)} \\
(2)  &  \left[\begin{matrix} \left(\begin{matrix} 0 \\ 0 \end{matrix} \right),
      & \left(\cO_{1}^{2} \setminus \begin{matrix} 0 & 0 \end{matrix} \right)
    \end{matrix} \right]  &   \cO_1 \rtimes \cO_{1}^{*}  \\
(3)  &  \left[\cO_{1}^{2} \setminus \begin{matrix} \left(\begin{matrix} 0 \\ 0
     \end{matrix} \right), & \left(\begin{matrix} 0 & 0 \end{matrix} \right)
    \end{matrix} \right] &   \cO_1 \rtimes \cO_{1}^{*}   \\
(4)  &  \left[\cO_{1}^{2} \setminus \begin{matrix} \left(\begin{matrix} 0 \\ 0
     \end{matrix} \right), & \left(\begin{matrix} 0 & \cO_{1}^{*} \end{matrix}
      \right)
    \end{matrix} \right]  &  \cO_1   \\
(5)  &  \left[\cO_{1}^{2} \setminus \begin{matrix} \left(\begin{matrix} 0 \\ 0
     \end{matrix} \right), & \left(\begin{matrix} u^{*} & \cO_1 \end{matrix}
      \right)
    \end{matrix} \right], u^{*} \in \cO_{1}^{*} &  \cO_{1}^{*}\\
\hline  
\end{array}
\end{displaymath}
Collecting all the pieces we get the desired result. 
\end{pf}
From the above discussion one easily obtains representation zeta
function for the group $\glO{4}{2}$. 
\begin{remark}
The representation zeta function $R_{G}(\mathcal{D})$ for $\D = 1$ gives the 
number of conjugacy classes of $G$. From above, we can easily obtain the number
of conjugacy classes of groups $\glO{2}{2}$, $\glO{3}{2}$ and $\glO{4}{2}$. Number of conjugacy classes of $\glO{2}{2}$ and $\glO{3}{2}$ is already known, see Avni-Onn-Prasad-Vaserstein~\cite{MR2543507}.  
\end{remark}
\begin{table}
\caption{\label{Table3} Group $\glO{4}{2}$}
The elements $\rho$, $\sigma$, $\tau$, and $\omega$ are primitive elements of $\F_q$,
$\F_{q^2}$, $\F_{q^3}$ and $\F_{q^4}$ respectively, such that $\rho = \sigma^{q+1} = \tau^{q^2
  + q +1} = \omega^{q^3+q^2+q+1}$ and $\sigma = \tau^{q^2 + 1}$. 
\[
\scriptsize{
\begin{array}{|c|c|c|c|}
\hline
 Type & \mathrm{ Number\; of \; similarity } &
     \mathrm{Isomorphism\; type \; } & Index \\
   &   \mathrm{classes \; of \;given}   &  \mathrm{of \; centralizer}
   & [Z_{\glO{4}{1}}(A) : \glO{4}{1}] \\
A   &  \mathrm{type} \;\; (n_{A}) & Z_{\glO{4}{1}}(A)   & \\
\hline
\left(
\begin{array}{cccc}
\rho^a & 0 & 0 & 0 \\
0 & \rho^a & 0 & 0 \\
0 & 0 & \rho^{a} & 0 \\
0 & 0 & 0 & \rho^{a} 
\end{array}
\right) 
            &    q    &       \glO{4}{1} & 1  \\
&&&\\
\left(
\begin{array}{cccc}
\rho^a & 0 & 0 & 0 \\
1 & \rho^a & 0 & 0 \\
0 & 0 & \rho^{a} & 0 \\
0 & 0 & 0 & \rho^{a} 
\end{array}
\right)         &   q        &   G_{(2,1,1)} & \parbox[]{2 cm}{$(q^2 +1)(q^3-1)(q+1)$} \\

&&&\\
\left(
\begin{array}{cccc}
\rho^a & 0 & 0 & 0 \\
1 & \rho^a & 0 & 0 \\
0 & 0 & \rho^{a} & 0 \\
0 & 0 & 1 & \rho^{a} 
\end{array}
\right)    &   q   &    \G_{(2,2)} & \parbox[]{2
  cm}{$q(q^4-1)(q^3-1)$} \\
&&&\\
\left(
\begin{array}{cccc}
\rho^a & 0 & 0 & 0 \\
1 & \rho^a & 0 & 0 \\
0 & 1 & \rho^{a} & 0 \\
0 & 0 & 0 & \rho^{a} 
\end{array}
\right)    &   q    &   G_{(3,1)}  & \parbox[]{2 cm}{$q^2(q^4-1)(q^3-1)(q+1)$}  \\
&&&\\
 \left(
\begin{array}{cccc}
\rho^a & 0 & 0 & 0 \\
1 & \rho^a & 0 & 0 \\
0 & 1 & \rho^{a} & 0 \\
0 & 0 & 1 & \rho^{a} 
\end{array}
\right)       &  q  &   \cO_{4}^{*} & \parbox[]{2 cm}{$q^3(q^4-1)(q^3-1)(q^2-1)$}   \\
&&&\\
 \left(
\begin{array}{cccc}
\rho^a & 0 & 0 & 0 \\
0 & \rho^a & 0 & 0 \\
0 & 0 & \rho^{a} & 0 \\
0 & 0 & 0 & \rho^{b} 
\end{array}
\right)    &      q(q-1)    &   \glO{3}{1} \times \cO_{1}^{*} &
\parbox[]{2 cm}{$q^3(q+1)(q^2+1)$}  \\
&&&\\
\left(
\begin{array}{cccc}
\rho^a & 0 & 0 & 0 \\
1 & \rho^a & 0 & 0 \\
0 & 0 & \rho^{a} & 0 \\
0 & 0 & 0 & \rho^{b} 
\end{array}
\right)    &    q(q-1)    &    G_{(2,1)} \times \cO_{1}^{*} &
\parbox[]{2 cm}{$q^3(q^2+1)(q+1)^2(q^3-1)$} \\
 &&&\\
\left(
\begin{array}{cccc}
\rho^a & 0 & 0 & 0 \\
1 & \rho^a & 0 & 0 \\
0 & 1 & \rho^{a} & 0 \\
0 & 0 & 0 & \rho^{b} 
\end{array}
\right)    &   q(q-1)   &   \cO_{3}^{*} \times \cO_{1}^{*} & \parbox[]{2 cm}{$q^4(q^4-1)(q^3-1)(q+1)$} 
\\
&&&\\
\left(
\begin{array}{cccc}
\rho^a & 0 & 0 & 0 \\
0 & \rho^a & 0 & 0 \\
0 & 0 & \rho^{b} & 0 \\
0 & 0 & 0 & \rho^{b} 
\end{array}
\right)     &     \frac{1}{2}q(q-1)   &  \glO{2}{1} \times \glO{2}{1}
& \parbox[]{2 cm}{$q^4(q^2+1)(q^2+q+1)$} \\
&&&\\
\left(
\begin{array}{cccc}
\rho^a & 0 & 0 & 0 \\
1 & \rho^a & 0 & 0 \\
0 & 0 & \rho^{b} & 0 \\
0 & 0 & 0 & \rho^{b} 
\end{array}
\right)           &    q(q-1) &   \cO_{2}^{*} \times \glO{2}{1} &
\parbox[]{2 cm}{$q^4(q^2+q+1)(q^4-1)$}  \\
&&&\\
\left(
\begin{array}{cccc}
\rho^a & 0 & 0 & 0 \\
1 & \rho^a & 0 & 0 \\
0 & 0 & \rho^{b} & 0 \\
0 & 0 & 1 & \rho^{b} 
\end{array}
\right)   &     \frac{1}{2} q(q-1)   &  \cO_{2}^{*} \times \cO_{2}^{*}
& \parbox[]{2 cm}{$q^4(q+1)(q^4-1)(q^3-1)$}
\\
&&&\\
\left(
\begin{array}{cccc}
\rho^a & 0 & 0 & 0 \\
0 & \rho^a & 0 & 0 \\
0 & 0 & \rho^{b} & 0 \\
0 & 0 & 0 & \rho^{c} 
\end{array}
\right)     &   \frac{1}{2} q(q-1)(q-2)   &   \glO{2}{1} \times \cO_{1}^{*}
\times \cO_{1}^{*} & \parbox[]{2 cm}{$q^5(q+1)(q^2+1)(q^2+q+1)$}  \\
&&&\\
\left(
\begin{array}{cccc}
\rho^a & 0 & 0 & 0 \\
1 & \rho^a & 0 & 0 \\
0 & 0 & \rho^{b} & 0 \\
0 & 0 & 0 & \rho^{c} 
\end{array}
\right)  &  \frac{1}{2}q(q-1)(q-2)  &  \cO_{2}^{*} \times \cO_{1}^{*} \times
\cO_{1}^{*} & \parbox[]{2 cm}{$q^5(q^2+1)(q+1)^2(q^3-1)$}\\
&&&\\
\hline
\end{array}
}
\]
\end{table}
\begin{table}
\[
\scriptsize{
\begin{array}{|c|c|c|c|}
\hline
Type & \mathrm{ Number\; of \; similarity } &
     \mathrm{Isomorphism\; type \; } & Index \\
    &   \mathrm{classes \; of \;given}   &  \mathrm{of \; centralizer}
    & [ Z_{\glO{4}{1}}(A) : \glO{4}{1}]   \\
 A   &  \mathrm{type} \;\;(n_{A})                        &
 Z_{\glO{4}{1}}(A) & \\
\hline

\left(
\begin{array}{cccc}
\rho^a & 0 & 0 & 0 \\
0 & \rho^b  & 0 & 0 \\
0 & 0 & \rho^{c} & 0 \\
0 & 0 & 0 & \rho^{d} 
\end{array}
\right)     &  \parbox[]{2 cm}{$\frac{1}{24} q(q-1)(q-2)(q-3)$} & \parbox[]{2.3 cm}{$\cO_1^{*} \times
\cO_1^{*} \times \cO_1^{*} \times \cO_1^{*}$} & \parbox[]{2
  cm}{$q^6(q^3+q^2+q+1)(q+1)(q^2+q+1)$}\\
&&&\\
\left(
\begin{array}{cccc}
\rho^a & 0 & 0 & 0 \\
0 & \rho^a & 0 & 0 \\
0 & 0 & \rho^{b} & 0 \\
0 & 0 & 0 & \sigma^{bq} 
\end{array}
\right)     &   \frac{1}{2} q^2 (q-1) & \glO{2}{1} \times \F_{q^2}^* &
\parbox[]{2 cm}{$q^5(q^2+1)(q^3-1)$}\\
&&&\\
 \left(
\begin{array}{cccc}
\rho^a & 0 & 0 & 0 \\
1 & \rho^a & 0 & 0 \\
0 & 0 & \sigma^{b} & 0 \\
0 & 0 & 0 & \sigma^{bq} 
\end{array}
\right)   &  \frac{1}{2} q^2 (q-1)  &  \cO_{2}^* \times \F_{q^2}^* & \parbox[]{2 cm}{$q^5(q^3-1)(q^4-1)$}\\
&&&\\
\left(
\begin{array}{cccc}
\rho^a & 0 & 0 & 0 \\
0 & \rho^b & 0 & 0 \\
0 & 0 & \sigma^{c} & 0 \\
0 & 0 & 0 & \sigma^{cq} 
\end{array}
\right)   &  \frac{1}{4} q^2 (q-1)^2  &   \cO_{1}^* \times \cO_{1}^* \times
\F_{q^2}^* & \parbox[]{2 cm}{$q^6(q+1)(q^2+1)(q^3-1)$} \\ 
&&&\\
\left(
\begin{array}{cccc}
\sigma^a & 0 & 0 & 0 \\
0 & \sigma^{aq} & 0 & 0 \\
0 & 0 & \sigma^{a} & 0 \\
0 & 0 & 0 & \sigma^{aq} 
\end{array}
\right)    &   \frac{1}{2} q (q-1)    &  \GL_{2}(\F_{q^2}) &
\parbox[]{2 cm}{$q^4(q-1)(q^3-1)$}\\
&&&\\
\left(
\begin{array}{cccc}
\sigma^a & 0 & 0 & 0 \\
0 & \sigma^{aq} & 0 & 0 \\
1 & 0 & \sigma^{a} & 0 \\
0 & 1 & 0 & \sigma^{aq} 
\end{array}
\right)          &   \frac{1}{2} q(q-1)   &  \mathbf F_{q^2} \times \mathbf
F_{q^2}^{*}  & \parbox[]{2 cm}{$q^4(q^4-1)(q^3-1)(q-1)$}  \\
&&&\\
\left(
\begin{array}{cccc}
\sigma^a & 0 & 0 & 0 \\
0 & \sigma^{aq} & 0 & 0 \\
0 & 0 & \sigma^{b} & 0 \\
0 & 0 & 0 & \sigma^{bq} 
\end{array}
\right)   &   \frac{1}{8} q(q-1) (q^2 -q -2) &  \F_{q^2}^* \times
\F_{q^2}^* & \parbox[]{2 cm}{$q^6(q^2+1)(q^3-1)(q-1)$} \\
&&&\\
\left(
\begin{array}{cccc}
\rho^a & 0 & 0 & 0 \\
0 & \tau^b & 0 & 0 \\
0 & 0 & \tau^{bq} & 0 \\
0 & 0 & 0 & \tau^{bq^2} 
\end{array}
\right)   &   \frac{1}{3} q^2 (q^2 -1)   &  \cO_{1}^* \times
\F_{q^3}^* & \parbox[]{2 cm}{$q^6(q^4-1)(q^2-1)$}\\
&&&\\
\left(
\begin{array}{cccc}
\omega^a & 0 & 0 & 0 \\
0 & \omega^aq & 0 & 0 \\
0 & 0 & \omega^{aq^2} & 0 \\
0 & 0 & 0 & \omega^{aq^3} 
\end{array}
\right)   &  \frac{1}{4} q(q^3 -1)    &   \F_{q^4}^{*} &
\parbox[]{2 cm}{$q^6(q-1)(q^2-1)(q^3-1)$} \\
\hline  
\end{array}
}
\]
\end{table}

\end{document}